\documentclass{article}
\usepackage{amsmath,amssymb,amscd}
\newtheorem{Def}{Definition}[section]
\newtheorem{Thm}{Theorem}[section]
\newtheorem{Lem}{Lemma}[section]
\newtheorem{Rem}{Remark}[section]
\newtheorem{Claim}{Claim}[section]

\begin{document}
\title{Spectral spread and non-autonomous Hamiltonian diffeomorphisms}
\author{Yoshihiro Sugimoto \thanks{Research Institute for Mathematical Sciences, Kyoto University} \thanks{ e-mail:sugimoto@kurims.kyoto-u.ac.jp \ \ MSC:53D05,53D35,53D40}} 
\date{}
\maketitle

\begin{abstract}
For any symplectic manifold ${(M,\omega)}$, the set of Hamiltonian diffeomorphisms ${\textrm{Ham}^c(M,\omega)}$ forms a group and ${\textrm{Ham}^c(M,\omega)}$ contains an important subset ${\textrm{Aut}(M,\omega)}$ which consists of time one flows of autonomous(time-independent) Hamiltonian vector fields on $M$. One might expect that ${\textrm{Aut}(M,\omega)}$ is a very small subset of ${\textrm{Ham}^c(M,\omega)}$. In this paper, we estimate the size of the subset ${\textrm{Aut}(M,\omega)}$ in ${C^{\infty}}$-topology and Hofer's metric which was introduced by Hofer. Polterovich and Shelukhin proved that the complement ${\textrm{Ham}^c\backslash \textrm{Aut}(M,\omega)}$ is a dense subset of ${\textrm{Ham}^c(M,\omega)}$ in ${C^{\infty}}$-topology and Hofer's metric if ${(M,\omega)}$ is a closed symplectically aspherical manifold where Conley conjecture is established (\cite{PS}). In this paper, we generalize above theorem to general closed symplectic manifolds and general convex symplectic manifolds. So, we prove that the set of all non-autonomous Hamiltonian diffeomorphisms ${\textrm{Ham}^c\backslash \textrm{Aut}(M,\omega)}$ is a dense subset of ${\textrm{Ham}^c(M,\omega)}$ in ${C^{\infty}}$-topology and Hofer's metric if ${(M,\omega)}$ is a closed or convex symplectic manifold without relying on the solution of Conley conjecture.
\end{abstract}

\section{Introduction}
\subsection{Background}
For any closed manifold $M$, the set of diffeomorphisms ${\textrm{Diff}(M)}$ forms a group which is called diffeomorphism group. This is an infinite-dimensional Lie group (non-Banach). It is well known that any one-parameter subgroup of ${\textrm{Diff}(M)}$ is generated by a vector field. More precisely, for any one-parameter subgroup ${f:\mathbb{R}\rightarrow \textrm{Diff}(M)}$, there is a vector field ${X\in \Gamma(TM)}$ which satisfies
\begin{equation*}
f(t)=\exp(tX)
\end{equation*}
where ${\exp:\Gamma(TM)\rightarrow \textrm{Diff}(M)}$ is a time $1$ flow of a vector field. From the inverse function theorem, one might expect that there is an open neighborhood of the zero section ${\mathcal{U}\subset \Gamma(TM)}$ such that 
\begin{equation*}
\exp:\mathcal{U}\longrightarrow \textrm{Diff}(M)
\end{equation*}
gives a diffeomorphism onto an open neighborhood of ${\textrm{Id}}$, but this is far from true (\cite{F}, \cite{M} Warning 1.6). This is an example that inverse function theorem does not hold for non-Banach manifolds.  This implies that the set of diffeomorphisms which are generated by vector fields ${\{\exp(X)\}_{X\in \Gamma(TM)}}$ is very small in ${\textrm{Diff}(M)}$. 

For a symplectic manifold ${(M,\omega)}$, one important subgroup of ${\textrm{Diff}(M)}$ is the Hamiltonian diffeomorphism group ${\textrm{Ham}^c(M,\omega)}$. In this case, the set of Hamiltonian diffeomorphisms of the form ${\exp(X)}$ is autonomous Hamiltonian diffeomorphisms ${\textrm{Aut}(M,\omega)}$. So one might expect that the set of autonomous Hamiltonian diffeomorphisms is very small in the Hamiltonian diffeomorphism group. Hofer introduced the so-called Hofer's metric on ${\textrm{Ham}^c(M,\omega)}$. Polterovich and Shelukhin proved that the complement of ${\textrm{Aut}(M,\omega)}$ contains an open dense subset in Hofer's metric which is also dense in ${C^{\infty}}$-topology for closed symplectically aspherical manifolds(\cite{PS}). Their proof can not be adapted to general symplectic manifolds  because it was based on the solution of the Conley conjecture for closed symplectically aspherical manifolds (\cite{SZ, G, GG}). In this paper, we generalize the above theorem to general closed symplectic manifolds and general convex symplectic manifolds without relying on the solution of the Conley conjecture.

\subsection{Main results}
Let ${(M,\omega)}$ be a symplectic manifold and let $C_c^{\infty}(S^1\times M)$ be the set of smooth periodic Hamiltonian functions on $M$ such that their supports are compact.
\begin{gather*}
C_c^{\infty}(S^1\times M)=\{H:S^1\times M\longrightarrow \mathbb{R} \ | \ \textrm{supp}(H)\subset S^1\times M \ \textrm{is \ compact}\}
\end{gather*}
In this paper, we define the Hamiltonian vector field of $H$ as follows.
\begin{equation*}
\omega(X_H,\cdot)=-dH
\end{equation*}
We denote the time $t$ flow of the vector field $X_H$ by $\phi_H^t$ and time $1$ flow by ${\phi_H}$. We call $\phi_H$ a Hamiltonian diffeomorphism generated by $H$ and we denote the set of all such Hamiltonian diffeomorphisms by ${\textrm{Ham}^c(M,\omega)}$.
\begin{equation*}
\textrm{Ham}^c(M,\omega)=\{ \phi_H \ | \ H\in C_c^{\infty}(S^1\times M) \}
\end{equation*}
If ${H\in C_c^{\infty}(H)}$ (in other words, $H$ does not depend on ${S^1}$), $H$ is called autonomous Hamiltonian function and $\phi_H$ is called autonomous Hamiltonian diffeomorphism.
\begin{equation*}
\textrm{Aut}(M,\omega)=\{\phi_H \ | \ H\in C_c^{\infty}(M)\}
\end{equation*}
We want to know the size of ${\textrm{Aut}(M,\omega)}$ in ${\textrm{Ham}^c(M,\omega)}$. In this paper, we prove  ${\textrm{Aut}(M,\omega)}$ is small with respect to Hofer's metric. Hofer's norm of a Hamiltonian function is defined as follows.
\begin{equation*}
||H||=\int_0^1\max H_t-\min H_tdt
\end{equation*}
Then, Hofer's norm of $\phi\in \textrm{Ham}^c(M,\omega)$ is defined by
\begin{equation*}
||\phi||=\inf\{||H|| \ | \ \phi_H=\phi, H\in C_c^{\infty}(S^1\times M)\}
\end{equation*}
and Hofer's metric on ${\textrm{Ham}^c(M,\omega)}$ is defined by
\begin{equation*}
\rho(\phi, \psi)=||\phi \psi^{-1}||.
\end{equation*}
It is known that ${(\textrm{Ham}^c(M,\omega),\rho)}$ is a metric space. ${(M,\omega)}$ is called symplectically aspherical if it satisfies the following conditions for any ${u:S^2\rightarrow M}$.
\begin{itemize}
\item $\omega(u)=\int_{S^2}u^*\omega=0$ 
\item $c_1(u)=\int_{S^2}u^*c_1=0$ where $c_1$ is the first Chern class of ${(M,\omega)}$
 \end{itemize}
Polterovich and Shelukhin proved the following theorem.
\begin{Thm}[Polterovich-Shelukhin\cite{PS}]
Let ${(M,\omega)}$ be a closed symplectically aspherical manifold. Then ${\textrm{Ham}^c(M,\omega)}$ contains a subset $W$ which satisfies the following conditions.
\begin{enumerate}
\item ${W\cap \textrm{Aut}(M,\omega)=\emptyset}$
\item $W$ is $C^{\infty}$-dense
\item $W$ is open and dense in the topology induced by Hofer's metric
\end{enumerate}
\end{Thm}

We generalize this theorem into Theorem 1.2. We define convex symplectic manifolds.
\begin{Def}
Let ${(M,\omega)}$ be a symplectic manifold. ${(M,\omega)}$ is called convex if there is a sequence of codimension $0$ submanifolds ${\{M_n\}_{n\in \mathbb{N}}}$ such that
\begin{itemize}
\item $\partial M_n\neq \emptyset$
\item $M_{n-1}\subset M_n$
\item $M=\bigcup_{n\in \mathbb{N}}M_n$
\item For any $n$, there is a vector field $X_n$ defined near ${\partial M_n}$ which is outward pointing on ${\partial M_n}$ and satisfies ${\mathcal{L}_{X_n}\omega=\omega}$ (Liouville vector field).
\end{itemize}
holds. 
\end{Def}
\begin{Thm}
\begin{enumerate}
\item 
Let ${(M,\omega)}$ be a closed symplectic manifold. Then, the set of all non-autonomous Hamiltonian diffeomorphisms ${\textrm{Ham}^c\backslash \textrm{Aut}}$ is ${C^{\infty}}$-dense in ${\textrm{Ham}^c(M,\omega)}$. Moreover, ${\textrm{Ham}^c(M,\omega)}$ contains a subset $W$ which satisfies the following conditions.
\begin{enumerate}
\item ${W\cap \textrm{Aut}(M,\omega)=\emptyset}$
\item $W$ is $C^0$-dense (not ${C^{\infty}}$-dense)
\item $W$ is open and dense in the topology induced by Hofer's metric
\end{enumerate}
\item
Let ${(M,\omega)}$ be a convex symplectic manifold. Then ${\textrm{Ham}^c(M,\omega)}$ contains a subset $W$ which satisfies the following conditions.
\begin{enumerate}
\item ${W\cap \textrm{Aut}(M,\omega)=\emptyset}$
\item $W$ is $C^{\infty}$-dense
\item $W$ is open and dense in the topology induced by Hofer's metric
\end{enumerate}
\end{enumerate}
\end{Thm}

The aim of this paper is to prove Theorem 1.2. In the second section, we explain the definition of spectral spread which was introduced by Polterovich and Shelukhin (\cite{PS}) and we prove a lemma (see Lemma 2.1 below) in the full generality in the third section. Our proof of Lemma 2.1 simplifies that of Polterovich and Shelukhin.

Polterovich and Shelukhin proved Theorem 1.1 by using spectral spread and the solution of Conley conjecture for symplectically aspherical manifolds. More precisely, they used Conley conjecture as follows. Let ${H\in C^{\infty}(S^1\times M)}$ be any non-degenerate Hamiltonian function. The solution of  Conley conjecture for closed symplectically aspherical manifolds implies  that there are infinitely many integers ${k\in \mathbb{Z}_{\ge2}}$ for which there exists a contractible ${k}$-periodic orbit of the Hamiltonian vector field ${X_{H}}$
\begin{equation*}
x(t):\mathbb{R}/k\mathbb{Z}\longrightarrow M
\end{equation*}
which satisfies ${x(0)\neq x(1)}$.
For general closed symplectic manifolds, we can not argue in the same way because this fact does not hold in general. So we need new ideas. For general convex symplectic manifolds,  the proof is a bit more complicated than that in closed case because spectral spread is not well-defined globally. More precisely, spectral spread is only defined for each ${M_n}$ ($M=\bigcup_{n\in \mathbb{N}}M_n$) as follows.
\begin{equation*}
w_k^{(n)}:C_c^0(S^1\times M_n)\longrightarrow \mathbb{R} 
\end{equation*}
So we only have a sequence of spectral spreads ${\{w_k^{(n)}\}_{n\ge 1}}$. We give a proof of Theorem 1.2 in the fourth section.

\section*{Acknowledgement}
The author thanks his supervisor, Professor Kaoru Ono, for many useful comments, discussions and encouragement. The author is supported by JSPS Research Fellowship for Young Scientists No. 201601854. The author also thanks referees for fruitful suggestions, especially for pointing out an error in our manuscript.

\section{Floer homology and spectral spread}
\subsection{Closed case}
We recall the definition of Floer homology and spectral spread and their properties. Let ${(M,\omega)}$ be a closed symplectic manifold. Let ${H\in C^{\infty}(S^1\times M)}$ be a Hamiltonian function. $H$ is called non-degenerate if 
\begin{equation*}
d\phi_H:T_xM\rightarrow T_xM
\end{equation*}
does not have $1$ as an eigenvalue for all ${x\in \textrm{Fix}(\phi_H)}$. We denote the space of contractible loops by ${\mathcal{L}(M)}$.
\begin{equation*}
\mathcal{L}(M)=\{x:S^1\rightarrow M \ | \ x \textrm{ \ is \ contractible}\}
\end{equation*}
We define Novikov covering of ${\mathcal{L}(M)}$ as follows.
\begin{equation*}
\widetilde{\mathcal{L}}(M)=\{(u,x) \ | \ x\in \mathcal{L}(M), u:D^2\rightarrow M, \partial u=x\}/\sim
\end{equation*}
where $D^2$ is a two dimensional disk and equivalence relation ${\sim}$ is defined as follows.
\begin{gather*}
(u,x)\sim(v,y)\Longleftrightarrow \begin{cases} x=y  \\ \omega(u\sharp \overline{v})=0  \\ c_1(u\sharp \overline{v})=0
\end{cases}
\end{gather*}
Here ${\overline{v}}$ is a disc with opposite orientation on the domain (unit disk) and ${u\sharp \overline{v}}$ is the glued sphere.
We denote the equivalence class which is represented by ${(u,x)}$ by ${[u,x]}$.
We define an action functional ${A_H}$ on ${\widetilde{\mathcal{L}}(M)}$ as follows.
\begin{equation*}
A_H([u,x])=-\int_{D^2}u^*\omega+\int_0^1H(t,x(t))dt
\end{equation*}
Let $P(H)$ be the set of contractible periodic orbits of $X_H$.
\begin{equation*}
P(H)=\{x\in \mathcal{L}(M) \ | \ \dot{x}(t)=X_{H_t}(x(t))\}
\end{equation*}
Let $\widetilde{P}(H)$ be the Novikov covering of ${P(H)}$.
\begin{equation*}
\widetilde{P}(H)=\{[u,x]\in \widetilde{\mathcal{L}}(M) \ | \ x\in P(H)\}
\end{equation*}
We denote the set of $S^1$-dependent compatible almost complex structures on $M$ by ${\mathcal{J}(M)}$. Let $H$ be a non-degenerate Hamiltonian function and let ${J\in \mathcal{J}(M)}$ be an almost complex structure. We define the Floer chain complex of ${(H,J)}$ as follows.
\begin{equation*}
CF(H,J)=\bigg\{\sum_{z\in \widetilde{P}(H), a_z\in \mathbb{Q}}a_z\cdot z \ \bigg| \ \forall C\in \mathbb{R}, \sharp\{z\in \widetilde{P}(H) \ | \ a_z\neq 0, A_H(z)>C\}< \infty \bigg\}
\end{equation*}
We consider the following equation on the cylinder.
\begin{gather}
u:\mathbb{R}\times S^1\longrightarrow M  \notag \\
\partial_su(s,t)+J(t,u(s,t))(\partial_tu(s,t)-X_{H_t}(u(s,t)))=0  \tag{$*$}
\end{gather}
For ${z_-=[v_-,x_-]}$, $z_+=[v_+,x_+]\in \widetilde{P}(H)$, we consider the following space.
\begin{equation*}
\widetilde{\mathcal{M}}(z_-,z_+,H,J)=\Bigg\{u:\mathbb{R}\times S^1\rightarrow M \ \Bigg| \ \begin{matrix}u \textrm{ \ satisfies \ }(*) \\ \lim_{s\rightarrow \pm \infty}u(s,t)=x_\pm(t) \\ [v_-\sharp u,x_+]=[v_+,x_+]\end{matrix}\Bigg\}
\end{equation*}

The above $\widetilde{\mathcal{M}}(z_-,z_+,H,J)$ has a natural $\mathbb{R}$-action. For any $u\in \widetilde{\mathcal{M}}(z_-,z_+,H,J)$ and ${s_0\in \mathbb{R}}$, ${s_0\cdot u}$ is defined as follows.
\begin{equation*}
s_0\cdot u(s,t)=u(s+s_0,t)
\end{equation*}
We take the quotient of this ${\mathbb{R}}$-action.
\begin{equation*}
\mathcal{M}(z_-,z_+,H,J)=\widetilde{\mathcal{M}}(z_-,z_+,H,J)/\mathbb{R}
\end{equation*}
By counting the $0$ dimensional part, we can define the boundary operator on the Floer chain complex (\cite{FO, LT}).
\begin{gather*}
\partial:CF(H,J)\rightarrow CF(H,J) \\
\partial(z_-)=\sum_{z_+\in \widetilde{P}(H)}\sharp \mathcal{M}(z_-,z_+,H,J)\cdot z_+
\end{gather*} 
$\partial$ satisfies $\partial \circ \partial=0$. Floer homology ${HF(H,J)}$ is the homology of the chain complex ${(CF(H,J),\partial)}$. The action functional $A_H$ give a filtration on ${CF(H,J)}$ as follows.
\begin{equation*}
CF^{<a}(H,J)=\{\sum a_z\cdot z \in CF(H,J)\ | \ A_H(z)<a \}
\end{equation*}
${CF^{<a}(H,J)}$ is a subcomplex of ${CF(H,J)}$. In other words, 
\begin{equation*}
\partial(CF^{<a}(H,J))\subset CF^{<a}(H,J)
\end{equation*}
is satisfied (see also \cite{S} Lemma2.1).

For ${a,b\in \mathbb{R}}$, we define the quotient
\begin{equation*}
CF^{[a,b)}(H,J)=CF^{<b}(H,J)/CF^{<a}(H,J).
\end{equation*} 
We define ${HF^{<a}(H,J)}$ and ${HF^{[a,b)}(H,J)}$ as the homology of ${CF^{<a}(H,J)}$ and ${CF^{[a,b)}(H,J)}$. We use this filtered Floer homology to define spectral spread. 

For $k\ge 2, k\in \mathbb{Z}$, we have the ${\frac{1}{k}}$-rotation of $\mathcal{L}(M)$.
\begin{gather*}
R_k:\mathcal{L}(M)\rightarrow \mathcal{L}(M)  \\
R_k(x)(t)=x\Big(t+\frac{1}{k}\Big)
\end{gather*}
This action also induces an action on ${\widetilde{\mathcal{L}}(M)}$ as follows.
\begin{gather*}
R_k: \widetilde{\mathcal{L}}(M)\rightarrow  \widetilde{\mathcal{L}}(M) \\
[u,x]\rightarrow [R_k(u),R_k(x)]
\end{gather*}
where $R_k(u)$ is ${\frac{1}{k}}$-rotation of the disk ${D^2}$ composed with $u$. Let ${H\in C^{\infty}(S^1\times M)}$ be a Hamiltonian function such that 
\begin{equation*}
H^{(k)}(t,x)=kH(kt,x)
\end{equation*}
is non-degenerate. 

For $J\in \mathcal{J}(M)$, we define $J^{+\frac{1}{k}}\in \mathcal{J}(M)$ by
\begin{equation*}
J^{+\frac{1}{k}}(t,x)=J\Big(t+\frac{1}{k},x\Big).
\end{equation*}
Then, $R_k$ induces an isomorphism between two chain complexes as follows.
\begin{gather*}
R_k:CF(H^{(k)},J) \rightarrow CF(H^{(k)},J^{+\frac{1}{k}})  \\
R_k\bigg(\sum_{z\in \widetilde{P}(H^{(k)})}a_z\cdot z\bigg)= \sum_{z\in \widetilde{P}(H^{(k)})}a_z\cdot R_k(z)
\end{gather*}
$R_k$ also induces a homomorphism
\begin{equation*}
R_k:HF^{[a,b)}(H^{(k)},J)\rightarrow HF^{[a+d,b+d)}(H^{(k)},J^{+\frac{1}{k}}) 
\end{equation*}
for any ${a<b}$ and ${d>0}$. For any ${a<b}$ and ${d>0}$, we fix a homotopy ${(K_s,J'_s)}$ (${s\in \mathbb{R}}$) between ${(H^{(k)},J)}$ and ${(H^{(k)},J^{+\frac{1}{k}})}$ as follows.
\begin{itemize}
\item ${(K_s,J'_s)=(H^{(k)},J)}$   holds for ${s\le -T}$ (for some ${T>0}$).
\item ${(K_s,J'_s)=(H^{(k)},J^{+\frac{1}{k}})}$ holds for  ${s\ge T}$.
\item ${\int_{\mathbb{R}}\max_{(t,x)\in S^1\times M}\frac{\partial}{\partial s}K_s(t,x)ds<d}$
\end{itemize}
For any ${z_-=[v_-,x_-],z_+=[v_+,x_+]\in \widetilde{P}(H^{(k)})}$, we define a moduli space
\begin{equation*}
\mathcal{N}(z_-,z_+,K_s,J'_s)
\end{equation*}
as a space of maps ${u:\mathbb{R}\times S^1\rightarrow M}$ which satisfies the following equations.
\begin{gather*}
\partial_su(s,t)+J'_{s,t}(\partial_tu(s,t)-X_{K_{s},t}(u(s,t)))=0  \\
\lim_{s\rightarrow \pm \infty}u(s,t)=x_{\pm}(t) ,\ \ \ \ \ [v_-\sharp u,x_+]=[v_+,x_+]
\end{gather*}
\begin{Def}
We define a continuation homomorphism ${I_{a,b,d}}$ as follows. First, we define a chain map ${I_{a,b,d,(K_s,J'_s)}}$.
\begin{gather*}
I_{a,b,d,(K_s,J'_s)}:CF^{[a,b)}(H^{(k)},J)\rightarrow CF^{[a+d,b+d)}(H^{(k)},J^{+\frac{1}{k}})  \\
I_{a,b,d,(K_s,J'_s)}(z_-)=\sum_{z_+\in \widetilde{P}(H^{(k)})}\sharp \mathcal{N}(z_-,z_+,K_s,J'_s)\cdot z_+
\end{gather*}
The chain homotopy class of ${I_{a,b,d,(K_s,J'_s)}}$ does not depend on the choice of the homotopy ${(K_s,J'_s)}$. So it induces a unique map on their homologies. We denote this map by ${I_{a,b,d}}$.
\begin{gather*}
I_{a,b,d}:HF^{[a,b)}(H^{(k)},J)\rightarrow HF^{[a+d,b+d)}(H^{(k)},J^{+\frac{1}{k}})  \\
\end{gather*}
\end{Def}

We define $S_k$ as follows.
\begin{equation*}
S_k=R_k-I_{a,b,d}
\end{equation*}
\begin{Def}[spectral spread \cite{PS}]
We define the spectral spread of $H$ as follows. First, assume that ${H^{(k)}}$ is non-degenerate.
\begin{equation*}
w_k(H)=\sup \Bigg\{d\ge 0 \ \Bigg| \ \begin{matrix}S_k:HF^{[a,b)}(H^{(k)},J)\rightarrow HF^{[a+d,b+d)}(H^{(k)},J^{+\frac{1}{k}})  \\ \textrm{is \ not \ }0\textrm{ \ for \ some \ }a<b \end{matrix}\Bigg\}
\end{equation*}
Lemma 2.1(3) below implies that ${w_k(H)}$ is Lipschitz continuous with respect to Hofer's norm of Hamiltonian functions. So we extend ${w_k(H)}$ for any continuous ${H\in C^0(S^1\times M)}$ as follows. Let ${\{H_i\}_{i\ge 1}}$ be a sequence of smooth Hamiltonian functions such that ${H_i^{(k)}}$ are non-degenerate and ${||H-H_i||\rightarrow 0}$ holds.
\begin{equation*}
w_k(H)=\lim_{i\rightarrow \infty} w_k(H_i)
\end{equation*}
\end{Def}

\begin{Rem}
Polterovich and Shelukhin defined $w_{k,\alpha}(H)$ for any free homotopy class $\alpha$ on symplectically atoroidal manifold ${(M,\omega)}$. They used Floer homology of free homotopy class $\alpha$ to define $w_{k,\alpha}$. We use only Floer homology of contractible periodic orbits for the sake of simplicity. Our arguments can be used also for non-contractible case $w_{k,\alpha}$, but we omit it.
\end{Rem}
We need the following properties of ${w_k}$.
\begin{Lem}
$w_k(H)$ satisfies the following properties.
\begin{enumerate}
\item $0\le w_k(H) \le k||H||$
\item $w_k(H)=w_k(K)$ if $\phi_H=\phi_K$. In other words, $w_k$ induces a map
\begin{equation*}
w_k:\textrm{Ham}^c(M,\omega)\rightarrow \mathbb{R}_{\ge0}
\end{equation*}
\item $|w_k(H)-w_k(K)|\le k||H-K||$
\item $w_k(\phi)=0$ if $\phi \in \textrm{Aut}(M,\omega)$.
\end{enumerate}
\end{Lem}
\begin{Rem}
\begin{itemize}
\item (1) follows from (3) and (4).
\item Polterovich and Shelukhin proved Lemma 2.1 for symplectically aspherical manifolds.
\end{itemize}
\end{Rem}

We do not need to prove (3) because the proof in \cite{PS} can be adapted to the general case. It was a standard argument in Floer theory. So we only prove (2) and (4) in the next section.

\subsection{Convex case}
A convex symplectic manifold is a union of compact codimension $0$ submanifolds such that their boundaries are non-empty and these boundaries are of contact type.  In this subsection, we will define the spectral spread for each of these submanifolds. So, we treat a compact symplectic manifold ${(M,\omega)}$ such that ${\partial M\neq \emptyset}$ holds in this subsection.   

We call ${\partial M}$ of contact type if there exists a vector field near ${\partial M}$ such that 
\begin{itemize}
\item $X$ is outward pointing on ${\partial M}$
\item $\mathcal{L}_{X}\omega=\omega$
\end{itemize}
hold. In this subsection, we assume that ${\partial M}$ has contact type. Then ${\alpha=\iota_X\omega|_{\partial M}}$ is a contact form on ${\partial M}$. A neighborhood of $\partial M$ is identified with ${(1-\tau,1]\times \partial M}$ whose symplectic form on ${(r,y)\in (1-\tau,1]\times \partial M}$ is ${d(r\alpha)}$. Symplectic completion ${(\widehat{M},\widehat{\omega})}$ is defined as follows (\cite{V}).
\begin{itemize}
\item $\widehat{M}=M\bigcup_{\partial M}[1,\infty)\times \partial M$
\item \begin{equation*}
\widehat{\omega}=\begin{cases} \omega & \textrm{on \ }M \\ d(r\alpha) & \textrm{on \ }(r,y)\in [1,\infty)\times \partial M \end{cases}
\end{equation*}
\end{itemize}
An almost complex structure $J$ on $\widehat{M}$ is called of contact type if it satisfies the following conditions.
\begin{itemize}
\item $J$ preserves $Ker(r\alpha)\subset T(\{r\}\times \partial M)$ on $\{r\}\times \partial M$
\item Let $X$ be a Liouville vector field on ${[1,\infty)\times \partial M}$ (${X(r,y)=r\frac{\partial}{\partial r}}$) and let $R$ be a Reeb vector field on ${\{r\}\times \partial M}$ with respect to the contact form ${r\alpha}$. $J$ satisfies ${J(X)=R}$ and ${J(R)=-X}$.
\end{itemize}

Let $T>0$ be the smallest period of periodic Reeb orbits of ${(\partial M,\alpha)}$ and we fix ${0<\epsilon<\frac{1}{k}T}$. we define the following family of pairs of a Hamiltonian function and a $S^1$-dependent contact type almost complex structure.
\begin{equation*}
\mathcal{H}=\Bigg\{(H,J) \ \Bigg| \ \begin{matrix}J \ \textrm{is\ a} \ S^1\textrm{-dependent \ compatible\ almost \ complex \ structure\ of \ contact\ type} \\ H\in C^{\infty}(S^1\times \widehat{M}) \\ H(t,(r,y))=-\epsilon (r-1), \forall(r,y)\in [1,\infty)\times \partial M  \end{matrix}\Bigg\}
\end{equation*}
For such $(H,J)$ we have the following Claim.
\begin{Claim}[\cite{AS,V}]
Let ${(H,J)\in \mathcal{H}}$ be any element of ${\mathcal{H}}$ such that $H$ is non-degenerate and let ${u:\mathbb{R}\times S^1\rightarrow \widehat{M}}$ be any solution of 
\begin{equation*}
\partial_su(s,t)+J(\partial_tu(s,t)-X_{H^{(k)}}(u(s,t)))=0 .
\end{equation*}
Then ${u(s,t)\in M}$ holds.
\end{Claim}
This claim implies that we can define ${HF^{[a,b)}(H^{(k)},J)}$ for ${(H^{(k)},J)\in \mathcal{H}}$ if ${H^{(k)}}$ is non-degenerate. Let ${(H,J)}$ be an element of ${\mathcal{H}}$ such that ${H^{(k)}}$ is non-degenerate. Then we can define the spectral spread ${\widehat{w}_k(H)}$ as in the closed case as follows.
\begin{equation*}
\widehat{w}_k(H)=\sup \Bigg\{d\ge 0 \ \Bigg| \ \begin{matrix}S_k:HF^{[a,b)}(H^{(k)},J)\rightarrow HF^{[a+d,b+d)}(H^{(k)},J^{+\frac{1}{k}})  \\ \textrm{is \ not \ }0\textrm{ \ for \ some \ }a<b \end{matrix}\Bigg\}
\end{equation*}
This definition does not depend on the choice of contact type almost complex structure $J$. Let ${\mathcal{H}_{\epsilon}}$ be a set of continuous Hamiltonian functions as follows.
\begin{equation*}
\mathcal{H}_{\epsilon}=\Bigg\{H\in C^0(S^1\times \widehat{M}) \ \Bigg| \ \begin{matrix} H(t,(r,y))=-\epsilon (r-1) \\ \textrm{on} \ (r,y)\in [1,\infty)\times \partial M
\end{matrix} \Bigg\}
\end{equation*}

${\widehat{w}_k}$ also satisfies Lipschitz continuity as follows.
\begin{equation*}
|\widehat{w}_k(H)-\widehat{w}_k(K)|\le k||H-K||
\end{equation*}
So we can extend ${\widehat{w}_k}$ on ${\mathcal{H}_{\epsilon}}$.

In this paper, we define the space of compact supported Hamiltonian functions ${C^0_c(S^1\times M)}$ as follows.
\begin{equation*}
C_c^0(S^1\times M)=\{H\in C^0(S^1\times M) \ | \ \textrm{supp}(H)\subset S^1\times \textrm{Int}(M)\}
\end{equation*}
For such ${H\in C^0_c(S^1\times M)}$, we define the canonical extension ${H_{\epsilon}\in \mathcal{H}_{\epsilon}}$ as follows.
\begin{equation*}
H_{\epsilon}(t,x)=\begin{cases}   H(t,x) & x\in M  \\ -\epsilon (r-1) & x=(r,y)\in [1,\infty)\times \partial M\end{cases}
\end{equation*}
We extend spectral spread on ${(M,\omega)}$ as follows.
\begin{Def}
For any ${H\in C_c^0(S^1\times M)}$, we define ${w_k(H)}$ by
\begin{equation*}
w_k(H)=\widehat{w}_k(H_{\epsilon}) .
\end{equation*}
\end{Def}

\begin{Rem}
The above $w_k$ also satisfies Lemma 2.1. So $w_k$ induces a map from 
\begin{equation*}
\textrm{Ham}^c(M,\omega)=\{\phi_H \ | \ H\in C_c^{\infty}(S^1\times M)\}
\end{equation*}
to ${\mathbb{R}_{\ge 0}}$. In other words, we have the following map.
\begin{equation*}
w_k:\textrm{Ham}^c(M,\omega)\longrightarrow \mathbb{R}_{\ge 0}
\end{equation*}
\end{Rem}

\section{Proof of Lemma 2.1}
\subsection{Proof of Lemma 2.1(2)}
We first prove (2). See also Proposition 5.3 in \cite{U} and \cite{PS} where similar arguments are used. Proof of Lemma 2.1 in convex case (${\partial M\neq \emptyset}$) is the same as in the closed case. So we assume ${(M,\omega)}$ is a closed symplectic manifold in this section. We assume $H,K\in C^{\infty}(S^1\times M)$ are Hamiltonian functions such that 
\begin{itemize}
\item $\phi_H=\phi_K$
\item $H^{(k)},K^{(k)}$ are non-degenerate
\end{itemize}
hold. Let $L\in C^{\infty}(S^1\times M)$ be a Hamiltonian function such that 
\begin{itemize}
\item $\phi_L=id$
\item $H^{(k)}(t,x)=L\sharp K^{(k)}(t,x)=L(t,x)+K^{(k)}(t,(\phi_L^t)^{-1}x)$
\end{itemize}
hold. Then, the flow ${\phi_L^t}$ satisfies ${\phi_L^t=\phi_{H^{(k)}}^t(\phi_{K^{(k)}}^t)^{-1}}$.

We first compare ${\widetilde{P}(H^{(k)})}$ and ${\widetilde{P}(K^{(k)})}$. The loop ${\phi_L^t}$ in ${\textrm{Ham}^c(M,\omega)}$ acts on ${\mathcal{L}(M)}$ as follows.
\begin{gather*}
f:\mathcal{L}(M)\rightarrow \mathcal{L}(M)  \\
f(x)(t)=\phi_L^t(x(t))
\end{gather*}
We fix two points ${[u,x],[v,y]\in \widetilde{\mathcal{L}}(M)}$ such that ${y=f(x)}$ holds. Then, there is a unique covering transformation 
\begin{equation*}
\widetilde{f}:\widetilde{\mathcal{L}}(M)\rightarrow \widetilde{\mathcal{L}}(M)
\end{equation*}
such that $\widetilde{f}([u,x])=[v,y]$ holds and the following diagram is commutative.

$$
\begin{CD}
\widetilde{\mathcal{L}}(M) @>\widetilde{f}>>  \widetilde{\mathcal{L}}(M)  \\
@V\pi VV  @V\pi VV   \\
\mathcal{L}(M) @>f>> \mathcal{L}(M)
\end{CD}
$$
We can also see that the difference
\begin{equation*}
A_{K^{(k)}}(z)-A_{H^{(k)}}(\widetilde{f}(z))
\end{equation*}
does not depend on $z\in \widetilde{\mathcal{L}}(M)$. This follows from the following calculations. It suffices to prove that the differential
\begin{gather*}
T\widetilde{\mathcal{L}}(M)\longrightarrow \mathbb{R}   \\
X\mapsto D(A_{K^{(k)}}(z)-A_{H^{(k)}}(\widetilde{f}(z)))X
\end{gather*}
vanishes. We fix ${z=[u,x]\in \widetilde{\mathcal{L}}(M)}$. Then, the tangent space ${T_z\widetilde{\mathcal{L}}(M)}$ can be identified with the space of vector fields along ${x\in \mathcal{L}(M)}$. In other words, 
\begin{equation*}
T_z\widetilde{\mathcal{L}}(M)\cong x^*TM
\end{equation*}
holds. We fix a vector field ${X(t)\in x^*TM}$. Then the following equalities hold.
\begin{gather*}
D(A_{K^{(k)}}(z))X=-\int_{S^1}\omega(X(t),\dot{x}(t))dt+\int_{S^1}dK_t^{(k)}\cdot X(t)dt 
\end{gather*}
\begin{align*}
&D(A_{H^{(k)}}(\widetilde{f}(z)))X \\
&=-\int_{S^1}\omega(d\phi_L^t\cdot X(t),\frac{\partial}{\partial t}(\phi_L^t(x(t))))dt+\int_{S^1}dH_t^{(k)}(d\phi_L^t\cdot X(t))dt   \\
&=-\int_{S^1}\omega(d\phi_L^t\cdot X(t), X_{L_t}+d\phi_L^t \cdot \dot{x}(t))dt \\
&+\int_{S^1}(dL_t+dK_t^{(k)}\circ d(\phi_L^t)^{-1})\circ d\phi_L^t\cdot X(t)dt  \\
&=-\int_{S^1}\{\omega(X(t),\dot{x}(t))+dL_t\circ d\phi_L^t(X(t))\}dt \\
&+\int_{S^1}\{dL_t\circ d\phi_L^t(X(t))+dK_t^{(k)}\cdot X(t)\}dt  \\
&=D(A_{K^{(k)}}(z))X
\end{align*}
The above equalities imply that ${D(A_{K^{(k)}}(z)-A_{H^{(k)}}(\widetilde{f}(z)))}$ vanishes. So, ${A_{K^{(k)}}(z)-A_{H^{(k)}}(\widetilde{f}(z))}$ does not depend on ${z}$ and we denote this constant by ${C}$.

Next, we compare moduli spaces
\begin{equation*}
\widetilde{\mathcal{M}}(z_-,z_+,K^{(k)},J_0)
\end{equation*}
and 
\begin{equation*}
\widetilde{\mathcal{M}}(\widetilde{f}(z_-),\widetilde{f}(z_+),H^{(k)},J_1)
\end{equation*}
for any $z_\pm \in \widetilde{P}(K^{(k)})$ and some $J_0,J_1\in \mathcal{J}(M)$. We fix ${J_1\in \mathcal{J}(M)}$. For ${u:\mathbb{R}\times S^1\rightarrow M}$, we define $v$ as follows.
\begin{gather*}
v:\mathbb{R}\times S^1\rightarrow M  \\
\phi_L^t(v(s,t))=u(s,t)
\end{gather*}
Then straightforward computation shows the following equation.
\begin{gather*}
\partial _su(s,t)+J_1(u(s,t))(\partial_tu(s,t)-X_{H^{(k)}}(u(s,t))) \\
=(\phi_L^t)_*(\partial_sv(s,t))+J_1(u(s,t))((\phi_L^t)_*(\partial_tv(s,t))-(
\phi_L^t)_*X_{K^{(k)}}(v(s,t)))
\end{gather*}
This implies that the equation 
\begin{equation*}
\partial_su+J_1(\partial_tu-X_{H{(k)}}(u))=0
\end{equation*}
and the equation
\begin{equation*}
\partial_sv+J_0(\partial_tv-X_{K^{(k)}}(v))=0
\end{equation*}
are equivalent for ${J_0=(\phi_L^t)_*^{-1}J_1(\phi_L^t)_*}$. We have the following ${1:1}$ correspondence.
\begin{gather*}
\widetilde{\mathcal{M}}(z_-,z_+,K^{(k)},J_0)\longrightarrow \widetilde{\mathcal{M}}(\widetilde{f}(z_-),\widetilde{f}(z_+),H^{(k)},J_1)  \\
v(s,t)\longrightarrow u(s,t)=\phi_L^t(v(s,t))
\end{gather*}
We can use $\widetilde{f}$ to identify ${CF^{[a,b)}(K^{(k)},J_0)}$ and ${CF^{[a+C,b+C)}(H^{(k)},J_1)}$. 

Next, we fix a homotopy ${(F_s,J'_s)}$ between ${(K^{(k)},J_0)}$ and ${(K^{(k)},J^{+\frac{1}{k}}_0)}$ and a homotopy ${(G_s,J''_s)}$ between ${(H^{(k)},J_0)}$ and ${(H^{(k)},J^{+\frac{1}{k}}_0)}$which satisfy the following conditions.
\begin{gather*}
G_s(t,x)=L\sharp F_s(t,s)  \\
J''_s=(\phi_L^t)_*^{-1}J'_s(\phi_L^t)_*
\end{gather*}

By using the similar calculation, we can see that there is the following ${1:1}$ correspondence.

\begin{gather*}
\mathcal{N}(z_-,z_+,F_s,J'_s)\longrightarrow \mathcal{N}(\widetilde{f}(z_-),\widetilde{f}(z_+),G_s,J''_s)  \\
v(s,t)\longrightarrow u(s,t)=\phi_L^t(v(s,t))
\end{gather*}

So we have the following commutative diagram.

$$
\begin{CD}
CF^{[a,b)}(K^{(k)},J_0) @>S_k>>  CF^{[a+d,b+d)}(K^{(k)},J^{+\frac{1}{k}}_0)  \\
@V\widetilde{f} VV  @V\widetilde{f}VV \\
CF^{[a+C,b+C)}(H^{(k)},J_1) @>S_k>>  CF^{[(a+C)+d,(b+C)+d)}(H^{(k)},J^{+\frac{1}{k}}_1)
\end{CD}
$$

This implies that ${w_k(H)=w_k(K)}$ holds. So we proved ${(2)}$.
\begin{flushright}     $\Box$ \end{flushright}

\subsection{Proof of Lemma 2.1(4)}
In this subsection, we prove (4). We fix a function 
\begin{equation*}
\rho:\mathbb{R}\longrightarrow \mathbb{R}
\end{equation*}
such that 
\begin{equation*}
\rho(s)=\begin{cases} 0 & s\le -1 \\ \frac{1}{k} & s\ge 1    \end{cases}
\end{equation*}
hold. We use $\rho$ to twist Floer equation as follows. For a cylinder ${u:\mathbb{R}\times S^1\rightarrow M}$, we define $\frac{1}{k}$-twist of $u$ by 
\begin{equation*}
v(s,t)=u(s,t+\rho(s)) .
\end{equation*}
We have the following equation.
\begin{gather*}
\partial_su(s,t)+J_t(\partial_tu(s,t)-X_{H^{(k)}}(u(s,t))) \\
=\partial_sv(s,t-\rho(s))-\rho'(s)\partial_tv(s,t-\rho(s)) \\ +J_t(v(s,t-\rho(s)))(\partial_tv(s,t-\rho(s))-X_{H^{(k)}}(v(s,t-\rho(s))))
\end{gather*}
This equation implies that the Floer equation
\begin{equation*}
\partial_su+J_t(\partial_tu-X_{H^{(k)}}(u))=0
\end{equation*}
is equivalent to
\begin{equation}
\partial_sv(s,t)-\rho'(s)\partial_tv(s,t)+J_{t+\rho(s)}(v)(\partial_tv(s,t)-X_{H_{t+\rho(s)}^{(k)}}(v(s,t)))=0 . \tag{A}
\end{equation}

We define the following moduli space for ${z_-,z_+\in \widetilde{P}(H^{(k)})}$.
\begin{equation*}
\mathcal{N}(z_-,z_+,\rho,H^{(k)},J)=\bigg\{v:\mathbb{R}\times S^1\rightarrow M \ \bigg| \ \begin{matrix} v \textrm{ \ satisfies \ }(\textrm{A})  \\ \lim_{s\rightarrow \pm \infty}v(s,t)=z_\pm, z_-\sharp u=z_+ \end{matrix} \bigg\}
\end{equation*}
Then, the above equivalence implies that there is an isomorphism between the following two moduli spaces.
\begin{equation*}
\widetilde{\mathcal{M}}(z_-,R_k^{-1}(z_+),H^{(k)},J)\cong \mathcal{N}(z_-,z_+,\rho,H^{(k)},J)
\end{equation*}
This ${\widetilde{\mathcal{M}}(\cdots)}$ is the space of connecting orbits which we defined in the second section. In particular, the $0$ dimensional part of ${\mathcal{N}(z_-,z_+,\rho,H^{(k)},J)}$ is compact and the following equality holds.
\begin{equation*}
\sharp \mathcal{N}(z_-,z_+,\rho,H^{(k)},J)=\begin{cases}  1  & z_+=R_k(z_-)  \\  0 & z_+\neq R_k(z_-)  \end{cases}
\end{equation*}
This equality implies that the following lemma holds.

\begin{Lem}
We define a chain map ${R_k^{'}}$ as follows.
\begin{gather*}
R_k^{'}:CF^{[a,b)}(H^{(k)},J)\longrightarrow CF^{[a,b)}(H^{(k)},J^{+\frac{1}{k}})  \\
z_-\longmapsto \sum \sharp \mathcal{N}(z_-,z_+,\rho,H^{(k)},J)\cdot z_+
\end{gather*}
Then, ${R_k^{'}=R_k}$ holds in chain level.
\end{Lem}

We deform the above equation (A) as follows.
\begin{equation}
\partial_sv(s,t)-\tau\cdot\rho'(s)\partial_tv(s,t)+J_{t+\rho(s)}(v)(\partial_tv(s,t)-X_{H_{t+\rho(s)}^{(k)}})=0 \tag{B}
\end{equation}
For ${z_-=[u_-,x_-],z_+=[u_+,x_+]\in \widetilde{P}(H^{(k)})}$, we define the following moduli space.
\begin{gather*}
\mathcal{N}(z_-,z_+,[0,1],\rho,H^{(k)},J) \\ =\Bigg\{(\tau,v)\in [0,1]\times C^{\infty}(\mathbb{R}\times S^1\rightarrow M) \ \Bigg| \ \begin{matrix}  v \textrm{ \ satisfies \ }(\textrm{B})  \\ \lim_{s\rightarrow \pm\infty}v(s,t)=x_{\pm}(t) \\  [u_-\sharp v,x_-]=[u_+,x_+]\end{matrix}\Bigg\}
\end{gather*}
The compactness of this moduli space follows from the following arguments. Let ${(\tau,v)}$ be an element of this moduli space. We twist $v$ as follows.
\begin{equation*}
u(s,t)=v(s,t-\tau \rho(s))
\end{equation*}
${(\tau,v)\in \mathcal{N}(z_-,z_+,[0,1],\rho,H^{(k)},J) }$ iff $u$ satisfies the following equation.
\begin{equation*}
\partial_su(s,t)+J_{t+(1-\tau)\rho(s)}(u)(\partial_tu(s,t)-X_{H_{t+(1-\tau)\rho(s)}^{(k)}}(u))=0
\tag{C}
\end{equation*}
Let ${R^{-\frac{1}{k}\tau}}$ be an obvious ${-\frac{1}{k}\tau}$ rotation as follows.
\begin{equation*}
R^{-\frac{1}{k}\tau}:\widetilde{P}(H^{(k)})\longrightarrow \widetilde{P}(H_{t-\frac{1}{k}\tau}^{(k)})
\end{equation*}
(Note that ${H_{t+\frac{1}{k}(1-\tau)}^{(k)}=H_{t-\frac{1}{k}\tau}^{(k)}}$ holds.) We define the following moduli space.
\begin{gather*}
\mathcal{M}(z_-,z_+,[0,1],\rho,H^{(k)},J) \\ =\Bigg\{(\tau,u)\in [0,1]\times C^{\infty}(\mathbb{R}\times S^1\rightarrow M) \ \Bigg| \ \begin{matrix}  u \textrm{ \ satisfies \ }(\textrm{C})  \\ \lim_{s\rightarrow \pm\infty}u(s,t)=x_{\pm}(t) \\  [u_-\sharp u,x_-]=R^{-\frac{1}{k}\tau}(z_+) \end{matrix}\Bigg\}
\end{gather*}
Then, there is an isomorphism as follows.
\begin{equation*}
\mathcal{N}(z_-,z_+,[0,1],\rho,H^{(k)},J)\cong \mathcal{M}(z_-,z_+,[0,1],\rho,H^{(k)},J)
\end{equation*}
The compactness of the $0$ dimensional part of  ${\mathcal{M}(z_-,z_+,[0,1],\rho,H^{(k)},J)}$ follows from the standard arguments in Floer theory. So ${\mathcal{N}(z_-,z_+,[0,1],\rho,H^{(k)},J)}$ is also compact.

Then, by counting the $0$ dimensional part of this moduli space, we construct $T$ as follows.
\begin{gather*}
T:CF(H^{(k)},J)\longrightarrow CF(H^{(k)},J^{+\frac{1}{k}})  \\
T(z_-)=\sum_{z_+\in \widetilde{P}(H^{(k)})}\sharp \mathcal{N}(z_-,z_+,[0,1],\rho,H^{(k)},J)\cdot z_+
\end{gather*}

Let ${(G_{s,t},J'_s)}$ be a homotopy between ${(H^{(k)},J)}$ and ${(H^{(k)},J^{+\frac{1}{k}})}$ as follows.
\begin{gather*}
G_{s,t}(x)=H^{(k)}_{t+\rho(s)} \\
J'_s(t,x)=J_{t+\rho(s)}(x)
\end{gather*}
Then Claim 3.1 below implies that for any ${d>k\max |\partial_tH|}$, $T$ induces a map
\begin{equation*}
T_{a,b,d}:CF^{[a,b)}(H^{(k)},J)\longrightarrow CF^{[a+d,b+d)}(H^{(k)},J^{+\frac{1}{k}})
\end{equation*}
and the homotopy ${(G_{s,t},J'_s)}$ determines  a chain map
\begin{equation*}
I_{a,b,d,(G_{s,t},J'_s)}:CF^{[a,b)}(H^{(k)},J)\longrightarrow CF^{[a+d,b+d)}(H^{(k)},J^{+\frac{1}{k}}) .
\end{equation*}
In other words, the homotopy ${(G_{s,t},J'_s)}$ can be used to define the map ${I_{a,b,d}}$ in Definition 2.1 and ${[I_{a,b,d,(G_{s,t},J'_s)}]=I_{a,b,d}}$ holds. 

Let ${\mathcal{N}^l(\cdots)}$ and ${\mathcal{M}^l(\cdots)}$ be connected components of moduli spaces ${\mathcal{N}}(\cdots)$ and ${\mathcal{M}(\cdots)}$ whose (virtual) dimensions are equal to $l$. Then, the boundary of ${\mathcal{N}^1(z_-,z_+,[0,1],\rho,H^{(k)},J)}$ can be written as follows.
\begin{gather*}
\partial \mathcal{N}^1(z_-,z_+,[0,1],\rho,H^{(k)},J)=\mathcal{N}^0(z_-,z_+,\rho,H^{(k)},J) \\
\bigsqcup \mathcal{N}^0(z_-,z_+,G_{s,t},J_s') \\
\bigsqcup \bigsqcup_{w\in \tilde{P}(H^{(k)})} \mathcal{M}^0(z_-,w,H^{(k)},J)\times \mathcal{N}^0(w,z_+,[0,1],\rho,H^{(k)},J)  \\
\bigsqcup \bigsqcup_{w\in \tilde{P}(H^{(k)})} \mathcal{N}^0(z_-,w,[0,1],\rho,H^{(k)},J)\times \mathcal{M}^0(w,z_+,H^{(k)},J^{+\frac{1}{k}})
\end{gather*}
(Here ${\mathcal{N}^0(z_-,z_+,G_{s,t},J_s')}$ is the moduli space which we used to define ${I_{a,b,d,(G_{s,t},J_s')}}$.)
This implies that 
\begin{equation*}
I_{a,b,d,(G_{s,t},J_s')}-R_k'=\partial \circ T_{a,b,d}+T_{a,b,d}\circ \partial
\end{equation*}
holds and ${T_{a,b,d}}$ is a chain homotopy between 
\begin{equation*}
R_k:CF^{[a,b)}(H^{(k)},J)\longrightarrow CF^{[a+d,b+d)}(H^{(k)},J^{+\frac{1}{k}})
\end{equation*}
and ${I_{a,b,d,(G_{s,t},J'_s)}}$.

\begin{Claim}
Assume that ${\mathcal{N}(z_-,z_+,[0,1],\rho,H^{(k)},J)\neq \emptyset}$ holds, then the following inequality holds.
\begin{equation*}
A_{H^{(k)}}(z_+)\le A_{H^{(k)}}(z_-)+k\max_{t\in [0,1]}|\partial_tH|
\end{equation*}
\end{Claim}
We fix $(\tau,v)\in \mathcal{N}(z_-,z_+,[0,1],\rho,H^{(k)},J)$. We define $u$ by
\begin{equation*}
u(s,t)=v(s,t-\tau\rho(s)) .
\end{equation*}
So ${(\tau, u)}$ is an element of ${\mathcal{M}(z_-,z_+,[0,1],\rho,H^{(k)},J)}$

We define ${w\in \widetilde{\mathcal{L}}(M)}$ by $z_-\sharp u$. The equation for ${(\tau,u)}$ implies that
\begin{gather*}
A_{H_{t-\frac{1}{k}\tau}^{(k)}}(w)\le A_{H^{(k)}}(z_-)+\int_{\mathbb{R}}\max_{t\in \mathbb{R}}|\frac{\partial}{\partial s}H^{(k)}_{t+(1-\tau)\rho(s)}|ds  \\
\le A_{H^{(k)}}(z_-)+(1-\tau)\int_{\mathbb{R}}\rho'(s)\max_{t\in \mathbb{R}}|\frac{\partial}{\partial t}H^{(k)}|ds \\ \le A_{H^{(k)}}(z_-)+k(1-\tau)\max|\partial_tH|
\end{gather*}
holds. On the other hand, we have
\begin{equation*}
A_{H^{(k)}_{t-\frac{1}{k}\tau}}(w)=A_{H^{(k)}}(z_+) .
\end{equation*}
So, we have
\begin{equation*}
A_{H^{(k)}}(z_+)\le A_{H^{(k)}}(z_-)+k\max_{t\in [0,1]}|\partial_tH| .
\end{equation*}
\begin{flushright}    $\Box$    \end{flushright}

We have proved that ${S_k=R_k-I_{a,b,d}}$ is zero and this implies that the following inequality holds.
\begin{equation*}
w_k(H)\le k\max|\partial_tH|
\end{equation*}
So, ${w_k(\phi)=0}$ holds for any ${\phi\in \textrm{Aut}(M,\omega)}$.
\begin{flushright}    $\Box$    \end{flushright}

\section{Proof of Theorem1.2}
\subsection{Closed case}
In this subsection, we assume that ${(M,\omega)}$ is a closed symplectic manifold. ${W=\bigcup_{k\ge 2}w_k^{-1}((0,\infty))}$ is a subset of ${\textrm{Ham}\backslash \textrm{Aut}(M,\omega)}$ and this subset is open in Hofer's metric. So it suffices to prove that ${W}$ is a  dense subset of  ${\textrm{Ham}^c(M,\omega)}$ in Hofer's metric and ${C^0}$-topology. We also prove that ${\textrm{Ham}\backslash \textrm{Aut}(M,\omega)}$ is ${C^{\infty}}$-dense in ${\textrm{Ham}(M,\omega)}$ independently. For these purposes, we use the following claim.

\begin{Claim}
If ${K^{(k)}}$ is non-degenerate and there is a contractible periodic orbit ${x\in P(K^{(k)})}$ such that ${x(0)\neq x(\frac{1}{k})}$ holds, then ${w_k(K)>0}$ and ${\phi_K\notin \textrm{Aut}(M,\omega)}$ hold.
\end{Claim}

Proof of this claim is straightforward. We fix ${z=[u,x]\in \widetilde{P}(K^{(k)})}$. Then ${R_k(z)\neq z}$ and ${A_{K^{(k)}}(z)=A_{K^{(k)}}(R_k(z))=a}$ hold for some ${a\in \mathbb{R}}$.We fix ${\epsilon>0}$ so that any non-trivial cylinder
\begin{gather*}
v:\mathbb{R}\times S^1\rightarrow M  \\
\partial_sv+J(\partial_tv-X_{K^{(k)}}(v))=0
\end{gather*}
satisfies 
\begin{equation*}
E(v)=\int_{\mathbb{R}\times S^1}\omega(\partial_sv,\partial_tv-X_{K^{(k)}}(v))dsdt>2\epsilon .
\end{equation*}
Then, ${\mathcal{M}(z,z_+,K^{(k)},J)\neq \emptyset}$ implies that ${A_{K^{(k)}}(z_+)<a-2\epsilon}$ holds and similarly, ${\mathcal{M}(z_-,R_k(z),K^{(k)},J^{+\frac{1}{k}})\neq \emptyset}$ implies ${A_{K^{(k)}}(z_-)>a+2\epsilon}$ holds. 

So ${z\in CF^{[a-2\epsilon,a+\epsilon)}(K^{(k)},J)}$ becomes a cycle. Next, we prove that ${S_k(z)\in HF^{[a-\epsilon, a+2\epsilon)}(K^{(k)},J^{+\frac{1}{k}})}$ is not zero. First, we fix a homotopy of ${S^1}$-dependent compatible almost complex structure ${J'_s}$ which satisfies the following conditions.
\begin{gather*}
J'_s(t,x)=\begin{cases}  J(t,x) & s<-R  \\ J^{+\frac{1}{k}}(t,x)=J(t+\frac{1}{k},x) & s> R \end{cases}
\end{gather*}
Let ${\{(G_{s,t},J'_s)\}_{m\in \mathbb{N}}}$ be a homotopy between ${(K^{(k)},J)}$  and ${(K^{(k)},J^{+\frac{1}{k}})}$ satisfying ${G_{s,t}(x)=K^{(k)}(x)}$ for all ${s\in \mathbb{R}}$, ${t\in S^1}$ and ${x\in M}$.
Then, ${A_{K^{(k)}}(z)=A_{K^{(k)}}(R_k(z))}$ and ${z\neq R_k(z)}$ implies that  
\begin{equation*}
\mathcal{N}(z,R_k(z),G_{s,t},J'_s)= \emptyset
\end{equation*}
holds. This implies that a representative of homology class ${S_k(z)\in HF^{[a-\epsilon, a+2\epsilon)}(K^{(k)},J^{+\frac{1}{k}})}$ can be written in the following form.
\begin{equation*}
R_k(z)-\sum_{w\in \widetilde{P}(K^{(k)})\setminus R_k(z), \ a-\epsilon \le A_{K^{(k)}}(w)<a+2\epsilon}\sharp \mathcal{N}(z,w,G_{s,t},J'_s)\cdot w
\end{equation*}
If ${z_-\in \widetilde{P}(K^{(k)})}$ satisfies ${A_{K^{(k)}}(z_-)<a+2\epsilon}$, ${\mathcal{M}(z_-,R_k(z),K^{(k)},J^{+\frac{1}{k}})=\emptyset}$. So the above chain does not become a boundary and ${S_k(z)\neq 0}$. This implies that ${w_k(K)\ge \epsilon}$ holds.
\begin{flushright}   $\Box$    \end{flushright}

We fix a Hamiltonian diffeomorphism ${\phi_H\in \bigcap_{k\ge 2}w_k^{-1}(0)}$. We construct a Hamiltonian function $K$ so that ${||H-K||}$ is arbitrary small.

We can perturb ${H}$ in ${C^{\infty}}$-topology so that ${\phi_H}$ and ${(\phi_H)^2}$ are non-degenerate. So without loss of generality, we can assume that ${\phi_H}$ and ${(\phi_H)^2}$ are non-degenerate. We fix ${\gamma(t)\in P(H)}$ and we denote ${\gamma(0)}$ by $p$. ${w_2(\phi_H)=0}$ and above claim imply that we can choose a small open neighborhood ${U_p}$ of $p$ so that the following condition holds.
\begin{equation*}
q\in U_p\backslash \{p\} \Longrightarrow \phi_H(q)\neq q, (\phi_H)^2(q)\neq q
\end{equation*}
We fix ${q\in U_p\backslash \{p\}}$ and a path ${l:[0,1]\rightarrow U_p}$ which satisfies the following conditions.
\begin{itemize}
\item ${l(0)=(\phi_H)^2(q)}$
\item ${l(1)=q}$
\item ${\phi_H(q)\notin \textrm{Im}(l)}$
\end{itemize}
We also choose a small open neighborhood ${U_l}$ of ${\textrm{Im}(l)}$ so that ${\phi_H(q)\notin U_l}$ holds and we choose a Hamiltonian function ${G}$ so that the following conditions hold.
\begin{itemize}
\item ${\textrm{supp}(G)\subset S^1\times U_l}$
\item ${(\phi_G\circ \phi_H)^2(q)=\phi_G\circ(\phi_H)^2(q)=q}$
\item ${d((\phi_G\circ \phi_H)^2)_{q}:T_qM\rightarrow T_qM}$ does not have $1$ as an eigenvalue.
\end{itemize}
Then, ${K'\stackrel{\mathrm{def.}}{=}G\sharp H}$ is a Hamiltonian function which satisfies the following conditions.
\begin{itemize}
\item ${\phi_{K'}(q)=\phi_H(q)\neq q}$
\item ${(\phi_{K'})^2(q)=\phi_G\circ (\phi_H)^2(q)=q}$
\item $q$ is a non-degenerate fixed point of ${(\phi_{K'})^2}$ and this orbit is contractible. 
\end{itemize}
In order to apply Claim 4.1, ${K'^{(2)}}$ have to be non-degenerate. Let ${K\in C^{\infty}(S^1\times M)}$ be a ${C^{\infty}}$-small perturbation of ${K'}$ so that 
\begin{equation*}
\textrm{supp}(K-K')\subset S^1\times M\backslash \{ \cup \{t\}\times \phi_{K'}^t(q)\bigcup \cup \{t\}\times \phi_{K'}^{1+t}(q) \}
\end{equation*}
holds and ${\phi_K}$ and ${(\phi_K)^2}$ are non-degenerate. Then ${\phi_K(q)\neq q}$ and ${(\phi_K)^2(q)=q}$ hold and ${\phi_K^{2t}(q)}$ (${t\in S^1}$) is contractible. Then Claim ${4.1}$ implies ${w_2(\phi_K)>0}$ holds. By making ${q\rightarrow p}$, we can make ${\phi_K}$ arbitrarily close to ${\phi_H}$ in ${C^0}$-topology and Hofer's metric. So, we proved that $W$ is dense in Hofer's metric and ${C^0}$-topology and open in Hofer's metric.

In the rest of this subsection, we prove that ${\textrm{Ham}\backslash \textrm{Aut}(M,\omega)}$ is ${C^{\infty}}$-dense in ${\textrm{Ham}(M,\omega)}$. We fix ${\phi \in \textrm{Aut}(M,\omega)}$. Without loss of generality, we assume that ${\phi}$ is strongly non-degenerate (In other words, ${\phi^k}$ is non-degenerate for every ${k\in \mathbb{Z}_{\ge 1}}$.). Non-degeneracy of ${\phi}$ implies that there is a Morse function ${f\in C^{\infty}(M,\mathbb{R})}$ such that ${\phi=\phi_f}$ holds. Let ${p\in M}$ be a critical point of ${f}$ so that 
\begin{equation*}
f(p)=\max_{x\in M}f(x)=C
\end{equation*}
holds. This implies that the Morse index of $p$ is equal to ${\textrm{dim}M=2n}$. Claim 4.1 and ${\phi \in \textrm{Aut}(M,\omega)}$ imply that 
\begin{equation*}
\textrm{Fix}(\phi^k)=\textrm{Fix}(\phi)
\end{equation*}
holds for any ${k\in \mathbb{Z}_{\ge 1}}$. Let $U$ be and open neighborhood of $p$ and let ${\psi}$ be a local chart (symplectic embedding) as follows.
\begin{equation*}
\psi:(U,\omega|_U)\longrightarrow (\mathbb{R}^{2n},\omega_0)
\end{equation*}
We assume that ${\psi(p)=(0,\cdots,0)}$ and ${\textrm{Fix}(\phi)\cap U=\{p\}}$ hold. Let ${(x_1,y_1,\cdots,x_n,y_n)}$ be local coordinates of ${\mathbb{R}^{2n}}$. Let $X$ be a vector field on $U$ as follows.
\begin{equation*}
X=\psi^*\Bigg(\frac{1}{2}\sum_i\bigg(x_i\frac{\partial}{\partial x_i}+y_i\frac{\partial}{\partial y_i}\bigg)\Bigg)
\end{equation*}
This $X$ satisfies ${\mathcal{L}_X\omega=\omega}$. We choose a sufficiently small ${a>0}$. Then,
\begin{equation*}
V(a)=\{x\in U \ | \ f(x)\ge C-a\}
\end{equation*}
is diffeomorphic to ${2n}$-disk and $X$ is outward pointing on ${\partial V(a)}$ (In other words, ${\partial V(a)}$ is of contact type with respect to the Liouville vector field $X$.). We prove the following claim.
\begin{Claim}
There is a ${T}$-periodic solution of the following equation for some ${T>0}$. 
\begin{equation*}
\begin{cases}
\gamma:\mathbb{R}\longrightarrow V(a)\backslash \{p\} \\
\frac{d\gamma}{dt}(t)=X_f(\gamma(t))
\end{cases}
\end{equation*}
So,
\begin{equation*}
\{x\in V(a)\backslash \{p\} \ | \ \phi_f^T(x)=x \ \textrm{for\ some\ }T>0\}\neq \emptyset
\end{equation*}
holds.
\end{Claim}
\begin{Rem}
The assumption ${\textrm{Fix}(\phi^k)=\textrm{Fix}(\phi)}$ and ${\textrm{Fix}(\phi)\cap U=\{p\}}$ implies that ${T\notin \mathbb{Q}}$ holds.
\end{Rem}
The proof of this claim relies on the computation of "local Floer homology" near ${p\in M}$. First, we assume that 
\begin{equation*}
\{x\in V(a)\backslash \{p\} \ | \ \phi_f^T(x)=x \ \textrm{for\ some\ }T>0\}= \emptyset
\end{equation*}
holds. We fix an almost complex structure $J$ on ${V(a)}$ which is of contact type on the boundary ${\partial V(a)}$. Then, Claim 2.1 implies that ${HF_*(f|_{V(a)}^{(l)}, J)}$ is well-defined for all ${l\in \mathbb{Z}_{\ge 1}}$ (Here, $*$ is the grading of Floer homology which comes from Conley-Zehnder index of periodic orbits.). We can compute this homology by using the assumption ${\textrm{Fix}(\phi_f^l)\cap U=\{p\}}$.
\begin{equation*}
HF_*(f|_{V(a)}^{(l)},J)\cong \begin{cases}  \mathbb{Q} & *=\textrm{Conley-Zehnder \ index \ of \ }p\in \textrm{Fix}(\phi_f^l)  \\ 0 & \textrm{others}
\end{cases}
\end{equation*}

Let ${\rho:\mathbb{R}\rightarrow \mathbb{R}}$ be a strictly increasing function which satisfies the following conditions. 
\begin{itemize}
\item ${\rho(s)=s}$ holds near ${s=C-a}$
\item ${\rho''(s)\le 0}$
\item ${\rho'(C)<1}$
\end{itemize}

Then, the Hamiltonian function ${g=\rho \circ f|_{V(a)}}$ on ${V(a)}$ satisfies ${\textrm{Fix}(\phi_g^l)=\{p\}}$ and 
\begin{equation*}
HF_*(g^{(l)},J)\cong \begin{cases}  \mathbb{Q} & *=\textrm{Conley-Zehnder \ index \ of \ }p\in \textrm{Fix}(\phi_g^l)  \\ 0 & \textrm{others}
\end{cases}
\end{equation*}
holds. The assumption ${\rho'(s)=1}$ near ${s=C-a}$ implies we can construct a homotopy ${(G_{s,t},J)}$ between ${(f|_{V(a)}^{(l)},J)}$ and ${(g^{(l)},J)}$ which satisfies the following conditions.
\begin{itemize}
\item ${G_{s,t}(x)=f|_{V(a)}^{(l)}(x)}$ holds on ${(s,t,x)\in (-\infty,-R]\times S^1\times V(a)}$ 
\item ${G_{s,t}(x)=g^{(l)}(x)}$ holds on  ${(s,t,x)\in [R,\infty)\times S^1\times V(a)}$
\item ${G_{s,t}(x)=f|_{V(a)}^{(l)}(x)=g^{(l)}(x)}$ hold on ${(s,t,x)\in \mathbb{R}\times S^1\times \mathcal{O}(\partial V(a))}$ (Here, ${\mathcal{O}(\partial V(a))}$ is an open neighborhood of ${\partial V(a)}$.)
\end{itemize} 
Then, by using Claim 2.1, we can see that a chain map
\begin{gather*}
CF_*(f|_{V(a)}^{(l)},J)\longrightarrow CF_*(g^{(l)},J)  \\
z_-\mapsto \sum_{z_+}\sharp \mathcal{N}(z_-,z_+,G_{s,t},J)\cdot z_+
\end{gather*}
is well defined and this induces an isomorphism between their homologies.
\begin{equation*}
HF_*(f|_{V(a)}^{(l)},J)\stackrel{\cong}{\longrightarrow} HF_*(g^{(l)},J)
\end{equation*}
This implies that the Conleys-Zehnder index of $p$ as a periodic orbit of ${\phi_f^l}$ is equal to the Conley-Zehnder index of $p$ as a periodic orbit of ${\phi_g^l}$. This is a contradiction because these two indices become different as ${l\in \mathbb{Z}_{\ge 1}}$ becomes large. So we proved the claim.

We choose ${q\in V(a)\backslash \{p\}}$ and ${T\in \mathbb{R}_{>0}}$ so that ${\phi_f^T(q)=q}$ holds. We also assume that ${T}$ is the smallest period (In other words, ${\phi_f(q)^t\neq q}$ holds for any ${0<t<T}$). Let ${\kappa>0}$ be a constant such that 
\begin{equation*}
\kappa T\in \mathbb{Q}\backslash \{\frac{1}{m} \ | \ m\in \mathbb{Z}_{\ge1}\}
\end{equation*}
hold. Let ${\rho:\mathbb{R}\rightarrow \mathbb{R}}$ be a strictly increasing function such that 
\begin{itemize}
\item ${\rho(s)=s}$ on ${s\le C-a}$
\item ${\rho'(f(q))=\frac{1}{\kappa}}$
\end{itemize}
hold. We define a Hamiltonian function ${f^{\rho}}$ by ${f^{\rho}=\rho \circ f}$. We also choose ${k\in \mathbb{Z}_{\ge 1}}$ so that 
\begin{equation*}
k(\kappa T)\in \mathbb{Z}
\end{equation*}
holds. Then, $1\notin \{m(\kappa T)\}_{m\in \mathbb{Z}_{>0}}$ implies that 
\begin{equation*}
\{q\}\in \textrm{Fix}(\phi_{f^{\rho}}^k)\backslash \textrm{Fix}(\phi_{f^{\rho}})
\end{equation*}
holds. We perturb ${f^{\rho}}$ as follows. Let ${F\in C^{\infty}(S^1\times M, \mathbb{R})}$ be a Hamiltonian function which satisfies the following conditions.
\begin{itemize}
\item ${\phi_F(\phi_{f^{\rho}}^m(q))=\phi_{f^{\rho}}^m(q)}$ holds for any ${m\in \mathbb{Z}}$
\item ${(F\sharp f_{\rho})^{(k)}}$ is non-degenerate
\end{itemize}
Then Claim 4.1 implies that ${w_k(\phi_{F\sharp f^{\rho}})>0}$ and ${\phi_{F\sharp f^{\rho}}\in \textrm{Ham}\backslash \textrm{Aut}(M,\omega)}$ hold. By making ${\kappa \rightarrow 1}$, we can make ${\rho(s)\rightarrow s}$ (in ${C^{\infty}}$-topology) and ${F\sharp f^{\rho}\rightarrow f}$ (in ${C^{\infty}}$-topology). So, we proved that ${\textrm{Ham}\backslash \textrm{Aut}(M,\omega)}$ is ${C^{\infty}}$-dense in ${\textrm{Ham}(M,\omega)}$.

\begin{flushright}   $\Box$    \end{flushright}
 
\subsection{Convex case}
In this subsection, we assume that ${(M,\omega)}$ is convex. We fix a sequence of codimension $0$ submanifolds ${\{M_m\}_{m\in \mathbb{N}}}$ such that 
\begin{itemize}
\item$M=\bigcup_{m\in \mathbb{N}}M_m$
\item $(M_m, \omega_m=\omega|_{M_m})$ is a symplectic manifold with a contact type boundary ${\partial M_m}$
\end{itemize}
holds. Then for any fixed ${k\in \mathbb{N}}_{\ge 2}$, we have a sequence of spectral spread ${w_k^{(m)}}$ as follows.
\begin{gather*}
w_k^{(m)}:\{H\in C^0_c(S^1\times M) \ | \ \textrm{supp}(H)\subset S^1\times \textrm{Int}(M_m) \}\longrightarrow \mathbb{R}
\end{gather*}
Let ${\delta>0}$ be a small positive real number so that ${k\delta}$ is smaller than the smallest period of periodic Reeb orbits on ${\partial M_m}$. 
Then, the definition of ${w_k^{(m)}(H)}$ is as follows.
\begin{equation*}
w_k^{(m)}(H)=\widehat{w}_k((H|_{M_m})_{\delta})
\end{equation*}
Above ${(H|_{M_m})_{\delta}\in C^0(S^1\times \widehat{M}_m)}$ is the canonical extension.

We also fix ${\phi \in \textrm{Ham}^c(M,\omega)}$. We will construct a subset ${U_{\phi}\subset \textrm{Ham}^c(M,\omega)}$ so that the following conditions hold.
\begin{itemize}
\item ${U_{\phi}\subset \textrm{Ham}^c\backslash \textrm{Aut}}$
\item ${U_{\phi}}$ is open with respect to Hofer's metric.
\item There is a sequence ${\phi_i\in U_{\phi}}$(${i=1,2,\cdots}$) such that ${\phi_i\rightarrow \phi}$ in ${C^{\infty}}$-topology and Hofer's metric.
\end{itemize}
Then ${W=\bigcup_{\phi \in \textrm{Ham}^c}U_{\phi}}$ is a ${C^{\infty}}$-dense and  an open dense subset of ${\textrm{Ham}^c}$ with respect to Hofer's metric. 

Let ${H\in C_c^{\infty}(S^1\times M)}$ be a Hamiltonian function which generates ${\phi}$ (${\phi=\phi_H}$). We choose a symplectic embedding 
\begin{equation*}
\iota:(B(3r),\omega_0)\hookrightarrow (M,\omega)
\end{equation*}
so that ${S^1\times \textrm{Im}(\iota)\cap \textrm{supp}(H)=\emptyset}$ holds. 

We fix a function ${\rho:[0,(3r)^2]\rightarrow [0,\infty)}$ and ${\epsilon >0}$ as follows.
\begin{gather*}
\rho(t)<\epsilon   \\ 
-\pi<\rho'(0)<0  \\
\textrm{supp}(\rho)\subset [0,(3r)^2) \\
\rho'(t)\le 0 
\\ \rho''(t)<0 \ \textrm{on} \ t< r^2  \\
\rho'|_{t\in [r^2,(2r)^2]}\equiv -\frac{4}{3}\pi  \ \ \Big(\Longrightarrow \min_{t\in [0,r^2]}\rho'(t)=-\frac{4}{3}\pi \Big)
\end{gather*}

\begin{Rem}
Let ${F(z)=\rho(|z|^2)}$ be a Hamiltonian function on ${B(3r)\subset \mathbb{C}^n}$. Then, $F$ generates the following Hamiltonian flow.
\begin{equation*}
\phi_F^t(z)=\exp(2\rho'(|z|^2)ti)\cdot z
\end{equation*}
\end{Rem}

We want to construct a sequence of Hamiltonian functions ${\{K_k\}_{k\in \mathbb{N}_{\ge 2}}}$ so that this sequence satisfies the following conditions.
\begin{itemize}
\item ${\textrm{supp}(H-K_{k})\subset S^1\times \iota(B(3r))}$
\item ${||H-K_k||<\frac{1}{k}\epsilon}$
\item There is a contractible periodic orbit ${x(t)\in P(K_k^{(k)})}$ which satisfies ${x(0)\neq x(\frac{1}{k})}$ and ${x(t)\in \iota(B(r))}$.
\item ${K_k^{(k)}}$ is non-degenerate on ${\iota(B(r))}$.
\end{itemize}

First, we define a sequence of Hamiltonian functions ${\{K'_k\}_{k \in \mathbb{N}_{\ge 2}}\subset C^{\infty}_c(S^1\times M)}$ as follows.
\begin{gather*}
K'_k(t,x)=H(t,x) \ \ (x\in M\backslash \iota(B(3r)))  \\
K'_k(t,\iota(z))=\frac{1}{k}\rho(|\iota(z)|^2) \ \ (z\in B(3r))
\end{gather*}
(We can define such ${K'_k}$ because the ball ${\iota(B(3r))}$ and the support of $H$ are disjoint.)
We define ${0<\tau<r}$ as follows.
\begin{gather*}
|\rho'(\tau^2)|=\pi
\end{gather*}

Such a $\tau$ is uniquely determined because we assumed ${\rho''(t)<0}$ on ${t< r^2}$. Then, ${x=\iota(0)}$ is the unique fixed point of ${\phi_{K_k'}}$ on ${\iota(B(r))}$ and ${x=\iota(z)}$ (${z\in B(r)}$) is a fixed point of ${\phi_{K'^{(k)}_k}}$ if and only if ${z=0}$ or ${|z|=\tau}$ hold. ${x=\iota(0)}$ is a non-degenerate fixed point, but ${x=\iota(z)}$ (${|z|=\tau}$) is a degenerate fixed point. We perturb ${K'_k}$ in a small neighborhood of the  sphere ${\iota(\partial B(\tau))}$ as follows. Let ${K_k}$ be a perturbed Hamiltonian function which satisfies the following conditions.
\begin{itemize}
\item ${||H-K_k||<\frac{1}{k}\epsilon}$
\item ${\iota(0)}$ is the unique fixed point of ${\phi_{K_k}}$ on ${\iota(B(r))}$
\item Fixed points of ${K_k^{(k)}}$ on ${\iota(B(r))}$ are non-degenerate.
\end{itemize}
Then, we have at least one contractible periodic orbit ${x\in P(K_k^{(k)})}$ such that ${x(0)\neq x(\frac{1}{k})}$ and ${x(t)\in \iota(B(r))}$ hold. This fact follows from the following arguments. ${\iota(0)}$ is the unique fixed point of ${\phi_{K_k}}$ on ${\iota(B(r))}$. So it suffices to prove that ${\phi_{K_k^{(k)}}}$ has at least two fixed points in ${\iota(B(r))}$. This follows from the following claim.

\begin{Claim}
Let ${T\in C^{\infty}(S^1\times B(r))}$ be a Hamiltonian function on ${(B(r),\omega_0)}$ which is defined as follows.
\begin{equation*}
T(t,z)=K_k(t,\iota(z))
\end{equation*}
Then, Floer homology of ${T^{(k)}}$ can be computed as follows. Here, ${*}$ is the grading of Floer homology which comes from Conley-Zehnder index of periodic orbits.
\begin{equation*}
HF_*(T^{(k)}, i)\cong \begin{cases} \mathbb{Q} & *=3n \\ 0 & *\neq 3n  \end{cases}
\end{equation*}
\end{Claim}

The Conley-Zehnder index of the periodic orbit ${\iota(0)}$ is $n$, so this claim implies that ${\phi_{K_k^{(k)}}}$ has at least two fixed point in ${\iota(B(r))}$.

This claim follows from the calculation of symplectic homology of balls in ${\mathbb{C}^{n}}$ (\cite{FHW}, Theorem 2). We fix ${\lambda \in \mathbb{R}_{>0}\setminus \{k\pi\}_{k\in \mathbb{N}}}$ and we fix a Hamiltonian function ${P\in C^{\infty}(S^1\times B(r))}$ such that 
\begin{equation*}
P(t,x)=-\lambda |x|^2+\textrm{constant}
\end{equation*}
holds on ${(t,x)\in S^1\times (B(r)\setminus X)}$ for some compact subset ${X\subset B(r)}$. Then, one can see that ${SH_*^{[-\infty, \lambda \pi^2)}(B(r))\cong HF_*(P)}$. In our case, ${\lambda=\frac{4}{3}\pi}$ and ${P=T^{(k)}}$. So, we can see that 
\begin{equation*}
HF_*(T^{(k)}, i)\cong SH_*^{[-\infty, \frac{4}{3}\pi r^2)}(B(r)) \cong \begin{cases} \mathbb{Q} & *=3n \\ 0 & *\neq 3n  \end{cases}
\end{equation*}

Next, we prove the theorem in the following two steps.
\begin{enumerate}
\item We prove that ${\phi_{K_k}\in \textrm{Ham}^c\backslash \textrm{Aut}}$ hold.
\item We construct an open neighborhood ${V_k}$ of ${\phi_{K_k}}$ (in Hofer's metric) so that ${V_k\subset \textrm{Ham}^c\backslash \textrm{Aut}}$ holds. 
\end{enumerate}
Then, we can define ${U_{\phi}}$ by ${U_{\phi}=\bigcup_{k\in \mathbb{N}_{\ge 2}}V_k}$ and ${\{\phi_{K_k}\in U_\phi\}}$ converges to ${\phi}$ in Hofer's metric and ${C^{\infty}}$-topology.\

We fix ${N\in \mathbb{N}}$ so that ${\textrm{supp}(K_k)\subset S^1\times \textrm{Int}(M_N)}$. For ${(1)}$, it suffices to prove that 
\begin{equation*}
\phi_{K_k}|_{M_m}\in \textrm{Ham}^c(M_m)\backslash \textrm{Aut}(M_m) \ \ (\forall m\ge N)
\end{equation*}
holds. We also fix ${m\ge N}$. So it suffices to prove that ${w^{(m)}_k(K_k)\neq 0}$ holds. Let ${\{z_j\}_{1\le j\le l}}$ be a subset of ${\widetilde{P}(K_k^{(k)})}$ which satisfies the following conditions. ${z_j}$ is written in the form ${[v,x]}$ such that 
\begin{gather*}
v:D^2\longrightarrow \iota(B(r))\subset M_m \\
x:S^1\longrightarrow \iota(B(r))\subset M_m
\end{gather*}
hold. We choose ${z=[v,x]\in \{z_j\}}$ which satisfies ${x(0)\neq x(\frac{1}{k})}$.

Let ${C_1>0}$ be a constant which is defined as follows.
\begin{equation*}
C_1=\min \{|A_{K_k^{(k)}}(z_i)-A_{K_k^{(k)}}(z_j)|\neq 0\}
\end{equation*}
Then
\begin{equation*}
|A_{K_k^{(k)}}(z)-A_{K_k^{(k)}}(z_j)|\ge C_1
\end{equation*}
or 
\begin{gather*}
A_{K_k^{(k)}}(z)=A_{K_k^{(k)}}(z_j)
\end{gather*}
hold. We denote ${A_{K_k^{(k)}}(z)}$ by ${a\in \mathbb{R}}$.

For a relatively compact connected open subset ${S\subset \mathbb{R}\times \mathbb{R}\slash 3\mathbb{Z}}$ such that ${\partial \overline{S}=\partial_1 S\bigsqcup \partial_2 S}$ (${\partial _1S}$ and ${\partial_2S}$ are not empty), we consider maps which satisfy the following conditions.
\begin{gather*}
f:\overline{S}\longrightarrow B(3r) \\
f(\partial_1 S)\subset \partial B(r), \  \ \ f(\partial_2 S)\subset \partial B(2r) \\
\partial_sf(s,t)+i\partial _tf(s,t)=0 \ \ \ \textrm{on} \ \ (s,t)\in S 
\end{gather*}
Then we define ${C_2}$ as follows. Let ${\omega_0}$ be a standard symplectic form on $\mathbb{C}^n$.
\begin{equation*}
C_2=\frac{1}{3}\inf \Big\{\int_Sf^*\omega_0  \ \Big| \ S \  \textrm{and} \ f \ \textrm{are as above}\Big\}
\end{equation*}

\begin{Rem}
This ${C_2}$ satisfies ${C_2>0}$. This follows from the famous monotonicity lemma (see for example, \cite{H} Theorem 1.3).
\end{Rem}

Let $C>0$ be a constant which satisfies ${C<\frac{1}{2}\min \{C_1,C_2\}}$. Let ${L\in C^{\infty}(S^1\times \widehat{M_m})}$ be a Hamiltonian function which satisfies the following conditions.
\begin{itemize}
\item ${L=K_k}$ on ${\iota(B(2r))}$
\item ${L(t,(r,y))=-\delta r}$ on ${(r,y)\in [1,\infty)\times \partial M_m}$
\item $L^{(k)}$ is non-degenerate.
\end{itemize}

The first condition implies that ${\{z_j\}_{1\le j\le l}}$ are also elements of ${\widetilde{P}(L^{(k)})}$. Next, we prove that  ${\widehat{w}_k^{(m)}(L)}\ge C$ holds. For this, it suffices to prove the following claim.

\begin{Claim}
Let ${z=[v,x]}$ be the element of ${\widetilde{P}(L^{(k)})}$ which we fixed above and let ${z_-,z_+}$ be any elements of ${\widetilde{P}(L^{(k)})}$. Then, the following two claims hold.
\begin{itemize}
\item $\mathcal{M}(z,z_+,L^{(k)},J)\neq \emptyset \Longrightarrow A_{L^{(k)}}(z_+)<a-2C$
\item $\mathcal{M}(z_-,R_k(z),L^{(k)},J^{+\frac{1}{k}})\neq \emptyset \Longrightarrow A_{L^{(k)}}(z_-)>a+2C  $
\end{itemize}
\end{Claim}

Then, by using the same arguments as in the proof of Claim ${4.1}$, we can see that $z\in \widetilde{P}(L^{(k)})$ becomes a cycle in ${CF^{[a-2C,a+C)}(L^{(k)},J)}$ and ${S_k(z)}$ is not zero in ${HF^{[a-C,a+2C)}(L^{(k)},J^{+\frac{1}{k}})}$.

Claim 4.3 follows from the following arguments. We fix an almost complex structure of contact type $J$ on ${\widehat{M}_m}$ which satisfies ${J|_{S^1\times \iota(B(3r))}=\iota^*i}$. Assume that ${\widetilde{\mathcal{M}}(z,z_+,L^{(k)},J)\neq \emptyset}$ holds. If ${z_+\in \{z_j\}_{1\le j \le l}}$ holds, then ${A_{L^{(k)}}(z_+)<a-2C}$ holds. If ${z_+\notin \{z_j\}_{1\le j \le l}}$ holds, ${u\in \widetilde{\mathcal{M}}(z,z_+)}$ satisfies the following conditions.
\begin{gather*}
u^{-1}(\iota(\partial (B(r))))\neq \emptyset , \ \ u^{-1}(\iota(\partial (B(2r))))\neq \emptyset
\end{gather*}
On ${u^{-1}(\iota(B(2r)\slash  B(r)))}$, a map ${g'(s,t)=\iota^{-1}(u(s,t))}$ satisfies the following equation.
\begin{equation*}
\partial_sg'(s,t)+i\{\partial_tg'(s,t)+\tfrac{8}{3}\pi i \cdot g'(s,t)\}=0
\end{equation*} 

In order to transform $g'$ into a holomorphic curve $f$ below, we consider the following three-fold covering.
\begin{equation*}
\pi:\mathbb{R}\times \mathbb{Z}\slash 3\mathbb{Z}\longrightarrow  \mathbb{R}\times S^1
\end{equation*}
Let ${S'\subset \mathbb{R}\times \mathbb{R}\slash 3\mathbb{Z}}$ be the preimage ${\pi^{-1}(u^{-1}(\iota(B(2r)\slash  B(r))))}$.
We define a map ${g:\overline{S'}\rightarrow B(3r)}$ by ${g(s,t)=\iota^{-1}(u(\pi(s,t)))}$. Then $g$ also satisfies the following Floer equation.
\begin{equation*}
\partial_sg(s,t)+i\{\partial_tg(s,t)+\tfrac{8}{3}\pi i \cdot g(s,t)\}=0
\end{equation*} 

Then, we define $f:\overline{S'}\rightarrow B(3r)$ as follows.
\begin{equation*}
f(s,t)=\exp(\tfrac{8}{3}\pi it)g(s,t)
\end{equation*}
This $f$ is a holomorphic curve as follows.
\begin{gather*}
\partial_sf(s,t)+i\partial_tf(s,t)=0 \\
f^{-1}(\partial B(r))\neq \emptyset, \ \ f^{-1}(\partial B(2r))\neq \emptyset
\end{gather*}
So we can choose a connected component ${S\subset S'}$ which satisfies the following conditions.
\begin{gather*}
\partial \overline{S}=\partial_1S\bigsqcup \partial_2S  \\
f(\partial_1 S)\subset \partial B(r), \  \ \ f(\partial_2 S)\subset \partial B(2r) \\
\partial_sf(s,t)+i\partial _tf(s,t)=0
\end{gather*}

The definition of $C_2>0$ implies that 
\begin{gather*}
A_{L^{(k)}}(z_+)=a-\int_{\mathbb{R}\times S^1} \omega(\partial_su,J\partial_su)dsdt \\
<a-\int_{u^{-1}(\iota(B(2r)\backslash B(r)))}\omega(\partial_su,J\partial_su)dsdt
\\=a-\frac{1}{3}\int_{S'}\omega_0(\partial_sg,i\partial_sg)dsdt   \\
=a-\frac{1}{3}\int_{S'}\omega_0(\partial_sf,i\partial_sf)dsdt
= a-\frac{1}{3}\int_{S'}f^*\omega_0 \\
\le a-\frac{1}{3}\int_{S}f^*\omega_0    \le a-C_2< a-2C
\end{gather*}
holds. 

This implies  that 
\begin{equation*}
\mathcal{M}(z,z_+,L^{(k)},J)\neq \emptyset \Longrightarrow A_{L^{(k)}}(z_+)<a-2C
\end{equation*}
holds. By using the same arguments, we can prove that the second claim holds. So, we proved  that ${\widehat{w}^{(m)}_k(L)\ge C}$ holds. By making ${||L-(K_k|_{M_m})_\delta||\rightarrow 0}$, we can see that 
\begin{equation*}
w^{(m)}_k(K_k)=\lim_{L\rightarrow (K_k|_{M_m})_\delta} \widehat{w}^{(m)}_k(L)\ge C>0
\end{equation*}
holds. So we proved ${(1)}$. For ${(2)}$, we define ${V_k\subset \textrm{Ham}^c(M,\omega)}$ as follows.
\begin{equation*}
V_k=\Big\{\psi \in \textrm{Ham}^c(M,\omega) \ \Big| \ ||\psi-\phi_{K_k}||<\frac{C}{k}\Big\}
\end{equation*}
So it suffices to prove the following claim.
\begin{Claim}
${V_k\subset \textrm{Ham}^c\backslash \textrm{Aut}}$ holds.
\end{Claim}
Recall that above constant ${C>0}$ does not depend on the choice of ${m\ge N}$. So we proved that ${w_k^{(m')}(K_k)}\ge C$ holds for any ${m'\ge N}$.

We fix ${\psi=\phi_G \in V_k}$ and ${m'\ge N}$ so that ${\textrm{supp}(G)\subset M_{m'}}$. Then,
\begin{equation*}
w^{(m')}_k(G)\ge w^{(m')}_k(K_k)-k\cdot \frac{C}{k}>C-C=0
\end{equation*}
holds. So ${\psi=\phi_G \notin \textrm{Aut}(M,\omega)}$ holds and we proved the theorem.
\begin{flushright}   $\Box$  \end{flushright}


\begin{thebibliography}{9}
\bibitem{AS}
M. Abouzaid, P. Seidel.
\textit{An open string analogue of Viterbo functoriality.}
Geometry $\&$ Topology 14 (2010), 627-718
\bibitem{FHW}
A. Floer, H. Hofer, K. Wysocki.
\textit{Applications of symplectic homology I.}
Mathematische Zeitschrift September 1994, Volume 217, Issue 1, pp 577-606
\bibitem{F}
C. Freifeld.
\textit{One-Parameter Subgroups Do Not Fill a Neighborhood of the Identity in an Infinite Dimensional Lie (Pseudo)-Group}.
Battelle Rencontres, 1967 ; Lectures in Mathematics and Physics (Benjamin, New York) p. 538-543
\bibitem{FO}
K. Fukaya, K. Ono.
\textit{Arnold conjecture and Gromov-Witten invariant.}
Topology Vol.38. No. 5, pp. 993-1048, 1999
\bibitem{G}
V.L. Ginzburg.
\textit{The Conley conjecture}.
Ann. of Math. 172 (2010), 1127-1180
\bibitem{GG}
V.L. Ginzburg, B.Z. G\"{u}rel.
\textit{On the generic existence of periodic orbits in Hamiltonian dynamics}.
J. Mod. Dyn. 3 (2009), 595-610
\bibitem{H}
C. Hummel.
\textit{Gromov's Compactness Theorem for Pseudo-holomorphic Curves}.
Progress in Mathematics Volume 151
\bibitem{LT}
G. Liu, G. Tian.
\textit{Floer homology and Arnold conjecture.}
J. Differential Geom. Volume 49, Number 1 (1998), 1-74
\bibitem{M}
J. Milnor.
\textit{Remarks on infinite-dimensional Lie groups}.
relativity, groups and topology Ⅱ, COURSE 10
\bibitem{PS}
L. Polterovich, E. Schelukhin.
\textit{Autonomous Hamiltonian flows, Hofer's geometry and persistence modules.}
Selecta Mathematica January 2016, Volume 22, Issue 1, pp227-296
\bibitem{SZ}
D. Salamon, E. Zehnder.
\textit{Morse theory for periodic solutions of Hamiltonian systems and the Maslov index}.
Comm. Pure Appl. Math. 45 (1992), 1303-1360
\bibitem{S}
M. Schwarz.
\textit{On the action spectrum for closed symplectically aspherical manifolds.}
Pacific Journal of Mathematics Vol. 193, No. 2, 2000 419-461
\bibitem{U}
M. Usher.
\textit{Hofer's metrics and boundary depth.}
Ann. Sci. \'Ec. Norm. Sup\'er. (4) 46 (2013), no. 1, 57--128
\bibitem{V}
C. Viterbo.
\textit{Functors and computations in Floer homology with applications Part 1.}
Geometric and Functional Analysis (1999) Volume 9, Issue 5, 985-1033
\end{thebibliography}
\end{document}